\documentclass[11pt]{article}

\usepackage{hyperref}

\usepackage{amssymb}

\usepackage{graphicx}

\usepackage{amssymb}
\usepackage{latexsym}
\usepackage{
enumerate}

\def\mathsf{\bf}

\def\d{\mathrm d}
\def\e{\mathrm e}

\def\R{\mathbb{R}}

\def\E{\mathrm E}

\def\text{\mbox}

{\catcode `\@=11 \global\let\AddToReset=\@addtoreset}
\AddToReset{equation}{section}

\newtheorem{theorem}{Theorem}[section]

\newtheorem{proposition}[theorem]{Proposition}

\newtheorem{remark}[theorem]{Remark}

\def\mathsf{\bf}

\def\R{\mathbb{R}}

\def\text{\mbox}

\def\e{\mathrm e}
\def\d{\mathrm{d}}

\begin{document}

\title{Covariance function of vector self-similar process \footnote{The joint work was partially supported by the research project MATPYL
of the F\'ed\'eration de Math\'ematiques des Pays de Loire.}}

\author{Fr\'ed\'eric Lavancier${}^{(1)}$, Anne Philippe${}^{(1)}$, Donatas Surgailis${}^{(2)}$\footnote{The author is partially supported by the Lithuanian State Science Foundation grant no. T-70/09. } }

\date{\today \\ \small \textbf{1.} Laboratoire de Math\'ematiques Jean Leray,
Universit\'e de Nantes, France\\
\texttt{\{Frederic.Lavancier,Anne.Philippe\}@univ-nantes.fr}  \\
\textbf{2.} Institute of Mathematics and Informatics, Vilnius,  Lithuania\\
 \texttt{sdonatas@ktl.mii.lt} \\
}

\maketitle

\begin{abstract}

The paper obtains the general form of the cross-covariance function of vector fractional Brownian motion with
correlated components having different self-simi\-la\-rity indices.

\end{abstract}
\begin{quote}

\end{quote}

\vskip2cm

\section{Introduction}

A $p-$variate stochastic process $X = \{X(t) = (X_1(t), \cdots, X_p(t)), t\in \R\} $ is said {\it operator  self-similar} (os-s) if there exists a $p\times p$ matrix $H$ (called the  exponent of $X$) such that for any $\lambda >0$,
\begin{eqnarray} \label{oH}
X(\lambda t)&=_{\rm fdd}&\lambda^H X(t),
\end{eqnarray}
where $=_{\rm fdd}$ means equality of finite-dimensional distributions, and the $p\times p$ matrix $\lambda^H $ is defined
by the power series
$\lambda^H = \e^{H \log \lambda} = \sum_{k=0}^\infty  H^k (\log \lambda)^k/k!. $
Os-s processes were
studied in Laha and Rohatgi (1981), Hudson and Mason (1982), Maejima and Hudson (1984), Sato (1991) and other papers. A Gaussian
os-s process with stationary increments is called operator fractional Brownian motion (ofBm). For $p=1$ the class of ofBm
coincides with fundamental class of fractional Brownian motions (fBm) (see e.g. Samorodnitsky and Taqqu (1994)). Recall that
a fBm with exponent $H \in (0,1)$ can be alternatively defined as a stochastically continuous Gaussian process
$X = \{X(t), t\in \R\}$ with zero mean and covariance
\begin{equation} \label{covX}
\E X(s) X(t) = \frac{\sigma^2}{2}(|s|^{2H} + |t|^{2H} - |t-s|^{2H}), \qquad t, s \in \R,
\end{equation}
where $\sigma^2 = \E X^2(1)$. The form of covariance of general ofBm seems to be unknown and may be quite complicated. The structure of ofBm and stochastic integral representations were studied in Didier and Pipiras (2008).

A particular case of os-s processes corresponds to {\it diagonal} matrix $H = {\rm diag}(H_1, \cdots, H_p)$. In this case,
relation (\ref{oH}) becomes
\begin{eqnarray} \label{vH}
(X_1(\lambda t), \cdots, X_p(\lambda t))&=_{\rm fdd}&\big(\lambda^{H_1} X_1(t), \cdots, \lambda^{H_p} X_p(t)\big).
\end{eqnarray}
Below, a $p-$variate process  $X$ satisfying (\ref{vH}) for any $\lambda >0$ will be called {\it vector self-similar} (vs-s) and
a stochastically continuous Gaussian vs-s process with stationary increments (si) will be called a {\it vector fractional Brownian motion }
(vfBm).

 Vs-s processes (in particularly, vfBm) seem  to be most useful for applications and statistical analysis of multiple time series.  
 They arise as limits of normalized  partial sums  of multivariate long memory processes with discrete time, see  Marinucci and Robinson (2000), Davidson and de Jong (2000), Chung (2002), Davidson and Hashimadze (2008), Robinson (2008),  Lavancier {\it et al.} (2009). 
 In particularly, vs-s processes appear in  the OLS estimation in multiple linear regression model  (Chung, 2002), 
 the multiple local Whittle estimation (Robinson, 2008) and  the 
two-sample testing for comparison of long memory parameters  (Lavancier {\it et al.}, 2009).

 Note from (\ref{vH}) that each component $X_i = \{X_i(t), t\in \R\}, \ i=1, \cdots, p$ of vs-s process is a (scalar) self-similar process, the fact which is not true for general os-s processes.

The present paper obtains the general form of the (cross-)covariance function of vs-s si process $X$ with finite variance and
exponent $H = {\rm diag}(H_1, \cdots, H_p), \  0< H_i < 1$. According to Theorem \ref{th:autosim} below, under some regularity condition, for any $i, j =1, \cdots, p, i\ne j$
with $H_i + H_j \neq 1$, there exist $c_{ij}, c_{ji} \in \R$ such that
for any $s, t \in \R$
\begin{equation}
{\rm cov}(X_i(s),X_j(t)) \ = \frac{\sigma_i\sigma_j}{2} \left\{ c_{ij}(s) |s|^{H_i+ H_j} + c_{ji}(t)
      |t|^{H_i+H_j} - c_{ji}(t-s) |t-s|^{H_i + H_j}\right\}, \label{covs-s}
\end{equation}
  where $\sigma^2_i := {\rm var}(X_i(1))$ and
  \begin{eqnarray}
    c_{ij}(t)&:=&\cases{c_{ij}, &$t>0$, \cr
    c_{ji}, &$t<0$.\cr}
    \label{cij}
  \end{eqnarray}
A similar expression (involving additional logarithmic terms) for the covariance ${\rm cov}(X_i(s),$  $X_j(t))$ is obtained  in the case
$H_i + H_j = 1$. We prove Theorem \ref{th:autosim} by deriving from vs-s si property
a functional equation for the cross-covariance  of the type studied in Acz\'el (1966)
and Acz\'el and Hossz\'u (1965) and using the result of the last paper (Theorem \ref{azcel} below) about the uniqueness of this equation.
Section \ref{sec:stoch-integr-repr} discusses the existence of vfBm with covariance as in (\ref{covs-s}). We start with a double sided
stochastic integral representation similar to Didier and Pipiras (2008):
\begin{equation}\label{XW}
X(t)=\int_{\R} \left\{ \left( (t-x)_+^{H-.5}- (-x)_+^{H-.5}\right) A_+ +
\left((t-x)_-^{H-.5} - (-x)_-^{H-.5}\right)A_- \right\} W({\d}x),
\end{equation}
where $ H-.5 := {\rm diag}(H_1-.5, \cdots, H_p - .5), x_+ := \max(x,0), x_- := \max(-x,0), A_+, A_- $ are real $p\times p$ matrices
and $W({\d}x) = (W_1({\d}x), \cdots, W_p({\d}x)) $ is a Gaussian white noise with zero mean, independent components and covariance
$\E W_i({\d}x) W_j({\d}x) = \delta_{ij} {\d}x $. According to Proposition 3.1, if $0<H_i< 1, H_i+H_j \ne 1, \ i,j=1, \cdots, p$ then
the cross-covariance of $X $ in (\ref{XW}) has the form as in (\ref{covs-s}) with
\begin{eqnarray} \label{tildeC1}
&&c_{ij}\ = \ 2 \tilde c_{ij} \phi_{ij}/\sigma_i \sigma_j, \quad
\phi_{ij} := B(H_i+ .5, H_j + .5)/\sin ((H_i+H_j)\pi),
\end{eqnarray}
where the matrix $\tilde C = (\tilde c_{ij}) $ is given by
\begin{eqnarray} \label{tildeC}
\tilde C&:=&\cos(H \pi) A_+ A_+^\ast + A_- A_-^\ast \cos(H \pi) \\
&-&\sin (H\pi) A_+ A_-^\ast \cos (H\pi) - \cos (H\pi) A_+ A_-^\ast \sin (H\pi). \nonumber
\end{eqnarray}
Here and below, $A^* $ denotes the transposed matrix,  $\sin (H \pi) := {\rm diag}(\sin (H_1 \pi),  $   $ \cdots, \sin (H_p \pi)), $  $
\cos (H \pi) := {\rm diag}(\cos (H_1 \pi), \cdots, \cos (H_p \pi))$.

\section{The form of the covariance function of vs-s process}

  Recall  that a random process $X = \{X(t), t\in \R \}$ has
  stationary increments (si)  if $\{X(t+T)$  $- X(T), t \in \R\}\
  =_{\rm fdd}\ \{X(t)- X(0), t \in \R\}$ for any $T \in
  \R$.

\begin{theorem}\label{th:autosim}
  Let $X = \left\{X(t), t\in \R \right\} $ be a 2nd order
  process with values in $\R^p$. Assume that
  $X$ has stationary increments, zero mean, $X(0) = 0$, and that $X$ is
  vector self-similar with exponent $ H = {\rm diag}(H_1, \cdots, H_p), \ 0< H_i < 1 \ (i=1, \cdots, p)$.

Moreover, assume also that for any $i, j = 1, \cdots, p$, the
  function $t \mapsto \E X_i(t) X_j(1)$ is continuously differentiable
  on $(0,1)\cup (1, \infty)$.  Let $\sigma_i^2>0$ denote the variance of $X_i(1), \ i= 1, \cdots, p$.

\smallskip

\noindent (i) If $i=j$, then for any $(s, t) \in \R^2$, we have
 \begin{equation}
    \E X_i(s)X_i(t) \ = \frac{\sigma_i^2}{2} \left\{|s|^{2H_{i}} +
      |t|^{2H_{i}} - |t-s|^{2H_{i}}\right\}. \label{cov1}
  \end{equation}

\noindent (ii) If $i \ne j$ and $H_i+H_j\not = 1$, then  there exists $c_{ij}, c_{ji} \in \R$ such that for any $(s, t) \in \R^2$,  (\ref{covs-s}) holds.

\noindent (iii) If $i\ne j$ and $H_i+H_j = 1$, then  there exists $d_{ij}, f_{ij} \in \R$ such that for any $(s, t )\in \R^2$, we have
    \begin{equation}\label{cov3}
      \E X_i(s)X_j(t) \ = \frac{\sigma_i\sigma_j}{2} \left\{ d_{ij}  (|s|+|t|-|s-t|)+ f_{ij} (t\log |t|-s\log |s| - (t-s) \log |t-s|)\right \}.
    \end{equation}

\noindent (iv) The matrix $R = (R_{ij})_{i,j=1,\cdots,p} $ is positive definite, where
$$
R_{ij}:=\cases{1, &$i=j$, \cr
c_{ij}+c_{ji}, &$i\ne j, \ H_i+ H_j \not = 1$, \cr
d_{ij}, &$i\ne j, \ H_i + H_j = 1$.\cr}
$$

\end{theorem}

\medskip

\noindent {\it Proof.} (i) Follows from the well-known characterization of covariance of (scalar-valued)
self-similar stationary increment process.

\smallskip

\noindent (ii)-(iii) From part (i), it follows that
$\E (X_i(t) - X_i(s))^2 = \sigma_i^2|t-s|^{2H_i}$ and hence $\{X_i(t), t \in \R\}$ is stochastically continuous on the real line, for any
$i=1, \cdots, p$.
Whence, it follows that $ \E X_i(s) X_j(t) $
is jointly continuous on $\R^2$ and vanishes for $t=0$ or $s=0$, for any $i,j=1, \cdots, p$.

Denote
$r(s,t) := \E X_i(s) X_j(t), \ H := \frac{1}{2}(H_i + H_j)$.  Let $\R_+ := \{u: u>0 \},\ \R_- := \{u: u<0\}$. From 
vs-s and si properties we obtain
\begin{eqnarray}
r(\lambda s, \lambda t)&=&\lambda^{2H} r(s,t), \label{cov-l} \\
r(s,t)&=&r(s+T,t+T)- r(s+T,T) - r(T,t+T) + r(T,T) \label{cov-T}
\end{eqnarray}
for any reals $s, t, T$ and any $\lambda >0$. Substituting $s=t=u, \lambda = 1/|u|$ into (\ref{cov-l}) one obtains
\begin{equation}
r(u,u) \ = \ 2\kappa_\pm |u|^{2H}, \quad u \in \R_{\pm}, \quad 2\kappa_\pm := r(\pm 1,\pm 1).  \label{cov-u}
\end{equation}
Substituting $s=t= -1, T = 1$ into (\ref{cov-T}) and using $r(0,0) = r(1,0) = r(0,1) = 0$ yields
$r(-1,-1) = r(1,1)$.
Next, substituting $T=u, s=t = v-u$ into (\ref{cov-T}) and using (\ref{cov-u}) with $\kappa := \kappa_+ = \kappa_-$ one obtains
\begin{eqnarray}
r(u,v)+ r(v,u)&=&r(v,v) + r(u,u) - r(v-u,v-u)\nonumber \\
&=&2\kappa(|v|^{2H} + |u|^{2H} - |u-v|^{2H}) \qquad (u,v \in \R).
\label{sym}
\end{eqnarray}
Next, let
$$
g_\pm(t)\ := \  r(\pm1,t)  \qquad { (t \in \R)}.
$$
Then (\ref{cov-l}) implies
\begin{equation}
r(s, t)\ = \  |s|^{2H} g_\pm({ \pm} t/s) \qquad (s\in \R_\pm,\,  { t \in \R}).
\label{SC1}
\end{equation}
Equation (\ref{cov-T}) with $s=1$ and (\ref{SC1}) yield, { for all $t\in\R$}
\begin{eqnarray}
g_+(t)&=&(T+1)^{2H} g_+\left(\frac{T+t}{T+1}\right)
- T^{2H} g_+\left(\frac{T+t}{T}\right)\nonumber \\
&-&(T+1)^{2H} g_+\left(\frac{T}{T+1}\right) + T^{2H} g_+(1) \qquad
 (T >0), \label{SC2}
\end{eqnarray}
and for all $t \in \R$,  equation (\ref{cov-T}) with $s=-1$ and (\ref{SC1}) yield
 \begin{eqnarray}
g_-(t)&=&(T-1)^{2H} g_+\left(\frac{T+t}{T-1}\right)
- T^{2H} g_+\left(\frac{T+t}{T}\right)\nonumber \\
&-&(T-1)^{2H} g_+\left(\frac{T}{T-1}\right) + T^{2H} g_+(1) \qquad
(T>1). \label{SC3}
\end{eqnarray}
We claim that the general solution of functional equations
(\ref{SC2})-(\ref{SC3}) has the form:

\begin{itemize}
\item When $H\not=1/2$,
\begin{eqnarray}\label{S1}
g_+(t)&=&\cases{c'+c'|t|^{2H}-c'|1-t|^{2H}, & \ \ $t<0$,\cr
c' + c |t|^{2H} - c'|1-t|^{2H}, & \ \ $0<t<1$, \cr
c' + c|t|^{2H} - c|1-t|^{2H}, & \ \ $t>1$, \cr}, \\
g_-(t)&=&\cases{c+c'|t|^{2H}-c'|t+1|^{2H}, &\ \ $t<-1$,\cr
c+c'|t|^{2H}-c|t+1|^{2H}, &\ \ $-1<t<0$,\cr
c + c|t|^{2H} - c|t+1|^{2H}, &\ \ $t>0$;\cr} \label{S2}
\end{eqnarray}

\item When $H=1/2$,
\begin{eqnarray}\label{S3}
g_+(t)&=&\cases{f\left(t\log |t|-(t-1)\log |t-1|\right), & \ $t<0$,\cr
d\left(1\land t\right) + f\left(t\log |t|-(t-1)\log |t-1|\right), &\  $t>0$, \cr}, \\
g_-(t)&=&\cases{d\left(1\land |t|\right) + f\left(t\log |t|-(t+1)\log |t+1|\right), & \ $t<0$,\cr
f\left(t\log |t|-(t+1)\log |t+1|\right), &\  $t>0$,\cr} \label{S4}
\end{eqnarray}

\end{itemize}
with some $c, c', d, f  \in \R$.

It follows from (\ref{S1})-(\ref{S4}) and (\ref{SC1}) that the
covariance $r(s,t) = \E X_i(s) X_j(t)$ for $(s, t) \in \R^2$ has the
form as in (\ref{covs-s}), (\ref{cov3}), with $c_{ij}= 2c/\sigma_i \sigma_j, c_{ji} = 2 c'/\sigma_i \sigma_j  $ in the case
$2H = H_i + H_j \ne 1 $, and
$d_{ij} = 2d/\sigma_i \sigma_j, \  f_{ij} = 2f/ \sigma_i \sigma_j$  in the case $2H = H_i + H_j  = 1 $.

\smallskip

To show the above claim, note, by direct verification,  that  (\ref{S1})-(\ref{S4}) solve equations (\ref{SC2})-(\ref{SC3}).
Therefore it suffices to show that (\ref{S1})-(\ref{S4}) is a unique solution of
(\ref{SC2})-(\ref{SC3}).

 Let $t>1$.  Differentiating (\ref{SC2}) with respect
to $t$ leads to
\begin{equation}
g'_+(t)\ =\  (T+1)^{2H-1}g'_+\Big(\frac{T+t}{T+1}\Big)- T^{2H-1}
g'_+\Big(\frac{T+t}{T}\Big). \label{dif}
\end{equation}
Let $x := t,\ y:= (T+t)/T$. Then $(x,y) \in (1,\infty)^2 $ and the
mapping $(t, T) \mapsto (x,y): (1, \infty)\times (0, \infty) \to
(1, \infty)^2 $ is a bijection. Equation (\ref{dif}) can be
rewritten as
\begin{equation}\label{dif1}
g'_+(F(x,y))\ =\
\Big(\frac{y-1}{x+y-1}\Big)^{2H-1}g'_+(x)+\Big(\frac{x}{x+y-1}\Big)^{2H-1}g'_+(y),
\end{equation}
where
\begin{equation}
F(x,y) \ := \  \frac{xy}{x+y-1}. \label{F}
\end{equation}

Equation (\ref{dif1}) belongs to the class of functional equations treated
in Acz\'el (1966) and Acz\'el and Hossz\'u (1965). For reader's convenience, we
present the result from Acz\'el (1966) which will be used below.

\bigskip

\begin{theorem} \label{azcel} {\rm (Acz\'el and Hossz\'u (1965))}  There exists at
    most one continuous function $f$ satisfying the functional
    equation
    \begin{equation}
      f(F(x,y)) \ = \  L(f(x), f(y), x, y)
    \end{equation}
    for all $x, y \in \langle A,B\rangle $ and the initial conditions
$$
f(a_1) = b_1, \quad f(a_2) = b_2 \qquad (a_1, a_2 \in \langle
A,B\rangle, \ a_1 \ne a_2 )
$$
if $F$ is continuous in $\langle A,B\rangle \times \langle A,B\rangle
$, $L(u,v,x,y)$ is strictly monotonic in $u$ or $v$ and $F$ is intern
(i.e., $F(x,y) \in (x,y) $ for all $x \ne y \in \langle A,B\rangle),$ where
$ \langle A,B\rangle $ is a closed, half-closed, open, finite or
infinite interval.
\end{theorem}
\bigskip

Unfortunately, equation (\ref{dif1}) does not satisfy the conditions of Theorem \ref{azcel}  (with  $f= g_+$ and $\langle A,B\rangle = (1, \infty)$), since
$F$ in (\ref{F}) is not intern. Following Acz\'el and Hossz\'u (1965), we first apply some transformations of (\ref{dif1}) so that
Theorem \ref{azcel} can be used.

Note,  taking $T = t/(t-1)$  in (\ref{dif}) one obtains
\begin{equation}
g'_+(t) \ = \  K(t) g'_+ \Big( \frac{t^2}{2t-1}\Big), \label{diff}
\end{equation}
where
\begin{equation}
K(t) \ := \ \frac{(2t-1)^{2H-1}}{t^{2H-1} + (t-1)^{2H-1}}.
\end{equation}
Let $\tilde F(x): =F(x,x)  =\frac{x^2}{2x-1}$. Then $\tilde F$ is
strictly increasing from $(1,\infty)$ onto $(1,\infty)$. Let
\begin{equation}
G(x,y) \  := \  \tilde F^{-1}(F(x,y)) \ = \  \frac{xy +
\sqrt{x(x-1)y(y-1)}}{x+y-1}, \label{G}
\end{equation}
with $\tilde F^{-1} (y) = y +\sqrt{y(y-1)}  $. Note $F(x,y) > 1, \
G(x,y) > 1 $ for $(x, y) \in (1, \infty)^2 $. From (\ref{diff})
with $t = G(x,y) $ one obtains
\begin{equation}
g'_+ (G(x,y)) \ = \  g'_+ (F(x,y)) K(G(x,y)). \label{dif2}
\end{equation}
Combining (\ref{dif1}) and (\ref{dif2}) one gets
\begin{eqnarray}\label{dif3}
g'_+(G(x,y))&=& L(g'_+(x),g'_+(y),x,y),
\end{eqnarray}
where
\begin{equation}
L(u,v,x,y)\ :=\  \left( \Big(\frac{y-1}{x+y-1}\Big)^{2H-1}u
+\Big(\frac{x}{x+y-1}\Big)^{2H-1}v \right)K(G(x,y)).
\end{equation}
The fact that $G$ in (\ref{G}) is intern follows from its definition and monotonicity in $x$ and $y$, implying
$$
x=G(x,x)\leq G(x,y) \leq G(y,y)=y
$$
for any $x\leq y, \ x,y \in (1,\infty)$. Since $L$ is monotonic in $u$ or $v$, so Theorem \ref{azcel} applies to functional equation
(\ref{dif3}) and therefore this equation has a unique continuous solution $g'_+ $ on the interval $(1, \infty)$, given boundary conditions $g'_+(a_i) = b_i \  (i=1,2), \  1 < a_2 < a_1 < \infty $.

\underline{Form of $g_+$ on $[1,+\infty)$} :
\begin{itemize}
\item { Assume $H\not=1/2$.} Let $a_1 := 2, b_1 := c H
  (2^{2H}-2)$. In view of (\ref{diff}), the other boundary condition
  can be defined by $a_2 = a_1^2/(2a_1-1) = 4/3, \ b_2 := g'_+(2)/K(2)
  $ and so equations (\ref{dif3}) and (\ref{dif}) have a unique
  solution for single boundary condition $g'_+(2) = cH
  (2^{2H}-2)$. (See also Acz\'el and Hossz\'u (1965, p.51).)  Since
  $g'_+(t) = 2H c(t^{2H-1} - (t-1)^{2H-1}) $ is a solution of
  (\ref{dif}) with this boundary condition, it follows that this
  solution is unique. Hence it also follows that
  \begin{equation}
    g_+(t) = c' + ct^{2H} - c(t-1)^{2H} \qquad t \in [1, \infty), \label{dif4}
  \end{equation}
  for some $c'\in \R$.

\item { When $H=1/2$. A particular solution of (\ref{dif}) when
    $t>1$ is $\log (t)-\log(t-1)$. For the same reason as above, the
    general solution of (\ref{dif}) is thus $d'(\log(t)-\log(t-1))$
    where $d' \in\R$. It follows that
    \begin{equation}\label{dif4half}
      g_+(t) = d + d'(t\log t - (t-1)\log (t-1)) \qquad t \in [1, \infty),
    \end{equation}
    for some $d, d' \in \R$.}

  \end{itemize}

\underline{Form of $g_+$ on $(0,1)$} :

{
Putting $t=1$ in  (\ref{SC2}) results in
$$
(T+1)^{2H} g_+\Big(\frac{T}{T+1}\Big) + T^{2H}g_+\Big(\frac{T+1}{T}\Big) \ = \
g_+(1)\big[(T+1)^{2H} + T^{2H} -1 \big].
$$
Whence, for $s:= T/(T+1) \in (0,1) $ and using (\ref{dif4}) one obtains

\begin{eqnarray}
g_+(s)&=&-s^{2H} g_+(1/s) + g_+(1)[1+ s^{2H} - (1-s)^{2H}]. \nonumber \\
\end{eqnarray}

This gives, for $s\in (0,1)$,
\begin{itemize}
\item when $H\not=1/2$,  $g_+(s)=c'+ cs^{2H} - c'(1-s)^{2H},$
\item when $H=1/2$,  $g_+(s)=d s + d' (s\log s - (s-1)\log (s-1)).$
\end{itemize}

Therefore, relations (\ref{S1}) and (\ref{S3}) have been proved when $t>0$. \\
}

\underline{Form of $g_+$ on $(-\infty,0)$} : {  This case follows from (\ref{SC2}) taking $T=-t.$}
\medskip

\underline{Form  of $g_-$ } :
The relations  (\ref{S2}) and (\ref{S4}) are deduced from (\ref{S1}) and (\ref{S3}) thanks to  relation (\ref{SC3}).

\smallskip

\noindent (iv) Follows from the fact that $R$ is
the covariance matrix of random vector $(X_1(1)/\sigma_1, $   $\cdots, X_p(1)/\sigma_p)$.

Theorem \ref{th:autosim} is proved.

\section{Stochastic integral representation of vfBm}
\label{sec:stoch-integr-repr}

In this section  we derive the covariance function of vfBm $X = \{X(t), t \in \R\}$ given by double-sided stochastic integral
representation in (\ref{XW}). Denote
\begin{eqnarray*}
\alpha^{++}_{ij} := \sum_{k=1}^p a^+_{ik} a^+_{jk}, \quad \alpha^{--}_{ij} := \sum_{k=1}^p a^-_{ik} a^-_{jk}, \quad
\alpha^{+-}_{ij} := \sum_{k=1}^p a^+_{ik} a^-_{jk}, \quad \alpha^{-+}_{ij} := \sum_{k=1}^p a^-_{ik} a^+_{jk},
\end{eqnarray*}
where $A_+ = \left(a^+_{ij}\right), A_- = \left(a^-_{ij}\right) $ are the $p\times p$ matrices in (\ref{XW}). Clearly,
$$
A_+ A_+^* = \left(\alpha^{++}_{ij}\right), \quad A_- A_-^* = \left(\alpha^{--}_{ij}\right), \quad A_+ A_-^* = \left(\alpha^{+-}_{ij}\right),
\quad A_- A_+^* = \left(\alpha^{-+}_{ij}\right).
$$
Note, each of the processes $X_{i} = \{X_i(t), t \in \R\}$ in (\ref{XW}) is a
well-defined fractional Brownian motion with index $H_i \in
(0,1)$; see e.g. Samorodnitsky and Taqqu (1994).

\begin{proposition} \label{lem:cov-intstoc} The covariance
of the
    process defined in (\ref{XW}) satisfies the following properties

  \begin{enumerate}[(i)]
\item   For any $i= 1, \cdots, p$  the variance of $X_i(1)$  is
  \begin{eqnarray}
    \label{C0}
    \sigma_i^2&=&\frac{B(H_i+.5, H_i + .5)}{\sin (H_i
        \pi)}  \left\{ \alpha^{++}_{ii} + \alpha^{--}_{ii}
      -  2\sin (H_i\pi) \alpha^{+-}_{ii} \right\}.
  \end{eqnarray}
\item  If $H_i+H_j \ne 1 $  then for
    any $s, t \in \R$, the cross-covariance $\E X_i(s) X_j(t) $ of the
    process in (\ref{XW}) is given by (\ref{covs-s}), with
    \begin{eqnarray}
     \frac{\sigma_i\sigma_j}{2} c_{ij}&:=&\frac{B(H_i+.5, H_j + .5)}{\sin ((H_i+H_j)
        \pi)} \nonumber \\
      && \hskip.2cm \times \left\{ \alpha^{++}_{ij} \cos(H_i \pi) + \alpha^{--}_{ij} \cos(H_j \pi)
      - \alpha^{+-}_{ij} \sin ((H_i+H_j)\pi)\right\}. \label{C}
    \end{eqnarray}

  \item  If $H_i+H_j  = 1 $  then   for any $s, t \in \R$, the
    cross-covariance $\E X_i(s) X_j(t) $ of the process in
    (\ref{XW}) is given by (\ref{cov3}), with
    \begin{eqnarray} \nonumber
      {\sigma_i\sigma_j} d_{ij}&:=&B(H_i+.5, H_j+.5) \\ && \hskip 1cm\times
      \bigg\{\frac{\sin(H_i\pi)+\sin(H_j\pi)}{2}(\alpha^{++}_{ij} + \alpha^{--}_{ij}) - \alpha^{+-}_{ij} - \alpha^{-+}_{ij} \bigg\}, \label{d-half} \\
       {\sigma_i\sigma_j} f_{ij}&:=&
(H_j-H_i)(\alpha^{++}_{ij} - \alpha^{--}_{ij}). \label{d-half2}
    \end{eqnarray}

  \end{enumerate}
\end{proposition}

\begin{remark} {\rm  Let $H_i + H_j \ne 1, \ i,j=1, \cdots, p$ and let $\tilde c_{ij}, \ \phi_{ij} $ be defined as
in (\ref{tildeC1}).

>From Proposition 3.1 (\ref{C0}), (\ref{C}) it follows that the matrix $\tilde C = (\tilde c_{ij}) $ satisfies  (\ref{tildeC}). In this
context,
a natural question arises  to find easily verifiable conditions on the matrices $\tilde C$ and $H$ such that there exist
matrices $A_+, A_- $ satisfying the quadratic matrix equation in (\ref{tildeC}).
In other words, for which $\tilde C$ and $H$ there exists
a vfBm $X$ with cross-covariance as in (\ref{covs-s})?

While the last question does not seem easy, it becomes much simpler if we restrict the class of vfBm's  $X$ in (\ref{XW})
to {\it causal} representations with $A_- = 0$. In this case, equation (\ref{tildeC}) becomes
$$
\tilde C = \cos(H) A_+ A_+^*.
$$
Clearly, the last factorization is possible if and only if the matrix $\cos(H)^{-1} \tilde C$ is symmetric and positive
definite.

}
\end{remark}

\bigskip

\noindent {\it Proof.} For all $s$, let
\begin{eqnarray*}
I^+_{ik}(s)&:=&\int_{\R} \left((s
-x)_+^{H_i-.5} - (-x)_+^{H_i-.5}\right) W_k(\d x),\\
I^-_{ik}(s)&:=&\int_{\R} \left((s
-x)_-^{H_i-.5} - (-x)_-^{H_i-.5}\right) W_k(\d x).
\end{eqnarray*}
Using the above notation,
$X_{i}(s)\ = \  \sum_{k=1}^p \left(a^+_{ik} I^+_{ik}(s) + a^-_{ik} I^-_{ik}(s)\right)$
and
\begin{eqnarray}\label{Isum}
\E X_i(s) X_j(t)&=&\alpha^{++}_{ij} \E I^+_{i1}(s) I^+_{j1}(t) + \alpha^{+-}_{ij} \E I^+_{i1}(s) I^-_{j1}(t) \\
&+&\alpha^{-+}_{ij} \E I^-_{i1}(s) I^+_{j1}(t) + \alpha^{--}_{ij} \E I^-_{i1}(s) I^-_{j1}(t).   \nonumber
\end{eqnarray}

Let $H_i+H_j \ne 1$. From Stoev and Taqqu (2006, Th. 4.1), taking there $a^+=1$, $a^-=0$, $H(s)=H_i$ and $H(t)=H_j$, we obtain
\begin{eqnarray*}
\E I^+_{i1}(s)I^+_{j1}(t)
&=&\psi_H \bigg[\cos\Big((H_j-H_i)\frac{\pi}{2}- \frac{(H_i + H_j)\pi}{2} \textrm{sign}(s)\Big)|s|^{H_i+H_j}\\
&&\ \  + \ \cos\Big((H_j-H_i)\frac{\pi}{2}+ \frac{(H_i + H_j)\pi}{2} \textrm{sign}(t)\Big)|t|^{H_i+H_j}\\
&&\ \  - \ \cos\Big((H_j-H_i)
\frac{\pi}{2}- \frac{(H_i + H_j)\pi}{2} \textrm{sign}(s-t)\Big)|s-t|^{H_i+H_j}\bigg],
\end{eqnarray*}
where
$$
\psi_H \ := \  \frac{\Gamma(H_i+.5)\Gamma(H_j+.5)\Gamma(2-H_i - H_j)}{\pi (H_i+H_j)(1-H_i-H_j)} \ = \
\frac{B(H_i+.5,H_j+.5)}{\sin ((H_i+H_j)\pi)}.
$$
Therefore,
\begin{eqnarray*}\E I^+_{i1}(s)I^+_{j1}(t)&=&\frac{B(H_i+.5,H_j+.5)}{\sin ((H_i+H_j)\pi)} \bigg[ b_{ij}(s) |s|^{H_i+H_j}+b_{ji}(t) |t|^{H_i+H_j}-b_{ij}(s-t) |s-t|^{H_i+H_j}\bigg],
\end{eqnarray*}
where
$$
b_{ij}(s)=\cases { \cos(H_i\pi), & if $s>0$ \cr
 \cos(H_j\pi), & if $s<0$.}
$$

Similarly, taking $a^+=0$, $a^-=1$, $H(s)=H_i$ and $H(t)=H_j$ in  Stoev and Taqqu (2006, Th. 4.1),  we obtain
\begin{eqnarray*}
\E I^-_{i1}(s)I^-_{j1}(t)&=&\frac{B(H_i+.5,H_j+.5)}{\sin ((H_i+H_j)\pi)}
\bigg[ b_{ji}(s) |s|^{H_i+H_j}+b_{ij}(t) |t|^{H_i+H_j}-b_{ji}(s-t) |s-t|^{H_i+H_j}\bigg] .\end{eqnarray*}

Finally,
\begin{eqnarray}
\E I^+_{i1}(s)I^-_{j1}(t)&=&
\Big( {\bf 1}_{\{s>t\}} \int_{t}^{s} (s-x)^{H_i-.5}(x-t)^{H_j-.5}{\rm d}x -
{\bf 1}_{\{s>0\}} \int_0^{s} (s-x)^{H_i-.5}x^{H_j-.5} \d x \nonumber\\ &&\hskip5cm  -  {\bf 1}_{\{t<0\}}
\int_{t}^0 (-x)^{H_i-.5} (x-t)^{H_j-.5} \d x \Big)\nonumber \\
&=&B(H_i+.5, H_j+.5)\Big( (s-t)_+^{H_i+H_j} - s_+^{H_i+H_j} - t_-^{H_i+H_j}\Big) \label{COV5}
\end{eqnarray}
and
\begin{eqnarray}
\E I^-_{i1}(s)I^+_{j1}(t)&=&B(H_i+.5, H_j+.5)\Big( (t-s)_+^{H_i+H_j} - t_+^{H_i+H_j} - s_-^{H_i+H_j}\Big). \label{COV6}
\end{eqnarray}
Substituting these formulas into (\ref{Isum}) we obtain (\ref{C}) and (\ref{C0}).

\smallskip

\noindent Next, let $H_i + H_j =1$. We get similarly from  Theorem 4.1 of Stoev and Taqqu (2006)
\begin{eqnarray*}
\E I^+_{i1}(s)I^+_{j1}(t)&=&  \frac{1}{\pi} B(H_i+.5, H_j+.5)
\bigg[\frac{\pi}{2}\sin(H_i\pi) (|s|+|t|-|s-t|)\\ && - \cos(H_i\pi) (s \log |s| - t \log |t|-(s-t) \log |s-t|)\bigg]
\end{eqnarray*}
and
\begin{eqnarray*}
\E I^-_{i1}(s)I^-_{j1}(t)&=&\frac{1}{\pi} B(H_i+.5, H_j+.5)\bigg[\frac{\pi}{2}\sin(H_i\pi) (|s|+|t|-|s-t|)\\ && + \cos(H_i\pi)
(s \log|s| - t \log |t|-(s-t) \log |s-t|)\bigg].
\end{eqnarray*}
Expressions (\ref{COV5}) and (\ref{COV6}) remain true when $H_i + H_j =1$ and
they can be rewritten  as
$$
\E I^+_{i1}(s)I^-_{j1}(t)\ = \  \E I^-_{i1}(s)I^+_{j1}(t)
\ = \
-\frac 12   B(H_i+.5, H_j+.5) (|s|+|t|-|s-t|).
 $$
Therefore, using (\ref{Isum}), we obtain (\ref{d-half}) and (\ref{d-half2}). Proposition 3.1 is proved.

\newpage

\section*{References}

\smallskip
\bigskip

\begin{description}
\itemsep -.04cm

\item{\sc Acz\'el, J.} (1966) {\em Lectures on functional equations and their applications.} Academic Press, New York.

\item{\sc Acz\'el, J. and Hossz\'u, M.} (1965) Further uniqueness theorems for functional equations.
{\em Acta. Math. Acad. Sci. Hung.} {\bf 16}, 51--55.

\item{\sc Chung, C.-F. } (2002) Sample means, sample autocovariances, and linear regression of stationary  multivariate long memory processes. {\em  Econometric Th.} {\bf 18}, 51--78.

\item{\sc Davidson, J. and de Jong, R.M. } (2000) The functional central limit theorem and weak convergence to stochastic
integrals. {\em  Econometric Th.} {\bf 16}, 643--666.

\item{\sc Davidson, J. and Hashimadze, N. } (2008) Alternative frequency and time domain  versions of fractional Brownian motion. {\em  Econometric Th.} {\bf 24}, 256--293.

\item{\sc Didier, G. and Pipiras,  V.} (2008) Integral representations of operator fractional Brownian motion. Preprint.

\item{\sc Hudson, W. and Mason, J.} (1982) Operator-self-similar processes in a finite-dimensional space.
{\em Trans. Amer. Math. Soc.} {\bf 273}, 281--297.

\item{\sc Laha, R.G. and Rohatgi, V.K.} (1981) Operator self-similar stochastic processes in $R^d$. {\em  Stochastic
Proces. Appl.} {\bf 12}, 73--84.

\item{\sc Lavancier, F., Philippe, A. and Surgailis, D.} (2009) A two-sample test for comparison of long memory parameters.
	arXiv:0907.1787v1 [math.ST]

\item{\sc Maejima, M. and Mason, J.} (1994) Operator-self-similar stable processes. {\em  Stochastic Proces.
Appl.} {\bf 54}, 139--163.

\item{\sc Marinucci, D. and Robinson, P.M.} (2000) Weak convergence of
 multivariate   fractional processes.
{\em Stochastic Process. Appl.} {\bf 86}, 103--120.

\item {\sc Robinson, P.M. } (2008) Multiple local Whittle estimation in stationary systems. {\em Ann. Statist.}
{\bf 36}, 2508--2530.

\item{\sc Samorodnitsky, G. and Taqqu, M.S.} (1994) {\em Sta\-ble Non-Gaussian Ran\-dom Pro\-ces\-ses.} \\
Chap\-man and Hall, New York.

\item{\sc Sato, K.} (1991) Self-similar processes with independent increments. {\em Probab. Th. Rel. F.} {\bf 89}, 285--300.

\item{\sc Stoev, S. and Taqqu, M.S.} (2006) How rich is the class of multifractional Brownian motions?
{\em Stochastic Process. Appl.}
{\bf 11}, 200--221.

\end{description}

\end{document}